# Robust mixture regression based on the skew *t* distribution

Fatma Zehra DOĞRU[1*] and Olcay ARSLAN[1]

[1]Ankara University, Faculty of Science, Department of Statistics, 06100 Ankara/Turkey
fzdogru@ankara.edu.tr, oarslan@ankara.edu.tr

**Abstract**

In this study, we propose a robust mixture regression procedure based on the skew *t* distribution to model heavy-tailed and/or skewed errors in a mixture regression setting. Using the scale mixture representation of the skew *t* distribution, we give an Expectation Maximization (EM) algorithm to compute the maximum likelihood (ML) estimates for the paramaters of interest. The performance of proposed estimators is demonstrated by a simulation study and a real data example.

## 1. Introduction

Mixture regression models are used to investigate the relationship between variables which come from some unknown latent groups. These models first introduced by Quandt (1972) and Quandt and Ramsey (1978) as switching regression models which are widely used in areas such as engineering, genetics, biology, econometrics and marketing. The parameter estimation of a mixture regression model is usally based on the normality assumption of the error terms. It is well-known that the estimators based on the normality assumption perform well when the error distribution is normal, but they are very sensitive to the departures (outliers, heavy-tailedness, skewness) from normality. To deal with the departures from normality robust mixture regression procedures have been proposed. Some of these works can be summarized as follows. Markatou (2000) and Shen et al. (2004) used a weight function to estimate the parameters robustly in the mixture regression models. Bashir and Carter (2012) used S-estimation method for the mixture linear regression model. Bai (2010) and Bai et al. (2012) proposed a robust estimation procedure based on M-regression estimation to estimate the parameters of the mixture linear regression model. Wei (2012) and Yao et al. (2014) explored the mixture regression model based on *t* distribution which is an extension of the mixture of *t* distribution studied by Peel and McLachlan (2000). Further, Zhang (2013) studied the robust mixture regression model using the Pearson Type VII distribution and Song et al. (2014) proposed a robust estimation procedure for the mixture regression models using the mixture of Laplace distribution as an error distribution. As it is pointed out by them, the robust mixture regression estimation procedure based on the Laplace distribution can be regarded as the application of the least absolute deviation (LAD) regression estimation to the mixture regression models. Liu and Lin (2014) proposed mixture regression model based on the skew normal distribution. Also, Pereira et al. (2012) studied performance of the estimates procedure for the mixtures of skew normal distribution.

In this paper, we propose a robust mixture regression procedure based on the skew *t* distribution to efficiently deal with heavy-tailedness and skewness in the mixture regression model setting. This is an extension of the mixture of skew *t* distribution proposed by Lin et al. (2007) to the mixture regression models. We will use the skew *t* distribution results from the scale mixture of the skew normal distribution introduced by Gupta et al. (2002), Gupta (2003) and Azzalini and Capitaino (2003). The scale mixture representation of the skew *t* distribution enables to easily implement an EM algorithm to obtain the ML estimates for the parameters of interest in the mixture regression model. One can see the works by Doğru and Arslan (2014) and Doğru (2015) on the mixture regression model based on the skew t distribution.

The paper is organized as follows. In Section 2, we give the basic definition of the mixture regression model. In Section 3, we present the robust mixture regression results based on the skew $t$ distribution. In Section 4 and 5, we give a simulation study and a real data example to compare the performances of the proposed estimation procedure with the other estimation procedures obtained from normal, $t$ (Yao et al. (2014)) and skew normal (Liu and Lin (2014)) distributions. The paper concludes with a conclusion section.

## 2. Mixture regression model

The model setting for a general mixture of linear regression models can be formulated as follows. Let $\boldsymbol{x}$ be a $p$-dimensional vector of explanatory variables, $y$ be the response variable and $Z$ be a latent class variable independent of $\boldsymbol{x}$. Suppose that given $Z = i$, the response variable $y$ depends on the explanatory variable $\boldsymbol{x}$ in a linear way

$$y = \boldsymbol{x}'\boldsymbol{\beta}_i + \epsilon_i, i = 1,2,\cdots,g, \qquad (1)$$

where, $\boldsymbol{\beta}_i = (\beta_{i1}, \beta_{i2}, \cdots, \beta_{ip})'$ is the unknown vector of regression parameters and $g$ is the number of components in mixture regression model. The random errors $\epsilon_i$ and $\boldsymbol{x}$ are assumed to be independent. In literature, it is often assumed that the random errors $\epsilon_i$'s have distributions from the location-scale family with zero means and $\sigma_i$ scale parameters. Suppose that $P(Z = i|\boldsymbol{x}) = w_i$, $i = 1,2,\cdots,g$, denote the mixing probabilities with $\sum_{i=1}^{g} w_i = 1$, then the conditional density function of $y$ given $\boldsymbol{x}$ can be of the form

$$f(y|\boldsymbol{x}, \boldsymbol{\Theta}) = \sum_{i=1}^{g} w_i f_i(y; \boldsymbol{x}'\boldsymbol{\beta}_i, \sigma_i), \qquad (2)$$

where, $f_i(y; \boldsymbol{x}'\boldsymbol{\beta}_i, \sigma_i)$ is the density function of the $ith$ component with some shape parameters (e.g. degrees of freedom for $t$ distribution) and $\boldsymbol{\Theta} = (w_1, \cdots, w_g, \boldsymbol{\beta}_1, \cdots, \boldsymbol{\beta}_g, \sigma_1, \cdots, \sigma_g)'$ is the unknown parameter vector. This model is called as a g-component mixture regression model.

The *ML* estimation method is used to estimate the unknown parameter vector $\boldsymbol{\Theta}$ in model (2). Let $\{(\boldsymbol{x_1}, y_1), (\boldsymbol{x_2}, y_2), \cdots, (\boldsymbol{x_n}, y_n)\}$ be a given sample. Then, the *ML* estimates is obtained by maximizing the following log-likelihood function with respect to $\boldsymbol{\Theta}$

$$\ell(\boldsymbol{\Theta}) = \sum_{j=1}^{n} \log\left[\sum_{i=1}^{g} w_i f_i(y_j; \boldsymbol{x}_j'\boldsymbol{\beta}_i, \sigma_i)\right]. \qquad (3)$$

However, it should be noted that the *ML* estimators cannot be explicitly obtained. The EM algorithm (Dempster et al., 1977) is used to find the *ML* estimates.

## 3. Robust mixture regression based on the skew *t* distribution

In this section, we will use the skew *t* distribution in order to model possible skewed and heavy-tailed errors in the mixture regression model. By doing so, we will obtain more robust estimators for the parameters of the mixture regression model. We will use the Azzalini type skew *t* distribution (Azzalini and Capitanio 2003, Gupta et al. 2002, Gupta 2003) with the following density function

$$f(\epsilon; \sigma^2, \lambda, \nu) = \frac{2}{\sigma} t_\nu(\eta) T_{\nu+1}\left(\lambda\eta\sqrt{\frac{\nu+1}{\eta^2+\nu}}\right), \eta = \frac{\epsilon}{\sigma}, \epsilon \in \mathbb{R}, \tag{4}$$

where, $\lambda \in \mathbb{R}$ is the skewness parameter, $t_\nu(\cdot)$ is the probability density function (pdf) of the $t$ distribution with $\nu \in (0, \infty)$ degrees of freedom and $T_{\nu+1}(\cdot)$ is the cumulative density function (cdf) of the $t$ distribution with $\nu + 1$ degrees of freedom.

In the mixture regression model (2), assume that the errors have a skew $t$ distribution with zero location, and $\sigma_i^2, \lambda_i$ and $\nu_i$ scale, skewness and degrees of freedom parameters, respectively. On the contrary to the symmetric case the mean $E(\epsilon_i) \neq 0$. For the skew t distribution $E(\epsilon_i) = \sigma_i \delta_{\lambda_i} \sqrt{\frac{\nu_i}{\pi}} \frac{\Gamma\left(\frac{\nu_i-1}{2}\right)}{\Gamma\left(\frac{\nu_i}{2}\right)}$ when $\nu_i > 1$, where $\delta_{\lambda_i} = \lambda_i/\sqrt{1+\lambda_i^2}$. Thus, $E(y_i) = x_j'\boldsymbol{\beta}_i + E(\epsilon_i)$, which only affects the intercept. Thus, when we estimate the intercept we will take into account this and correct $\hat{\beta}_0$ by using $\widehat{E(\epsilon_i)}$. In order to estimate the unknown parameters we should maximize the following log-likelihood function

$$\ell(\Theta) = \sum_{j=1}^{n} \log\left(\sum_{i=1}^{g} w_i f_i(y_j - x_j'\boldsymbol{\beta}_i; \sigma_i^2, \lambda_i, \nu_i)\right) \tag{5}$$

where, $\Theta = (w_1, \cdots, w_g, \boldsymbol{\beta}_1, \cdots, \boldsymbol{\beta}_g, \sigma_1^2, \cdots, \sigma_g^2, \lambda_1, \cdots, \lambda_g, \nu_1, \cdots, \nu_g)'$.

However, the maximizer of the above log-likelihood function cannot be explicitly obtained so that an EM-type algorithm should be used to estimate the unknown parameters $\Theta$. The EM algorithm can be implemented as follows. Let $Z_{ij}$ be the latent variables such that

$$Z_{ij} = \begin{cases} 1, & \text{if } j^{th} \text{ observation is from } i^{th} \text{ component} \\ 0, & \text{otherwise} \end{cases}$$

where, $j = 1, \cdots, n$ and $i, \cdots, g$. To simplify the EM algorithm we use the stochastic representation of a skew $t$ distributed random variable given by Azzalini and Capitanio (2003) (see Appendix for more details). This stochastic representation yields the following hierarchical formulation in terms of the conditional distributions

$$y_j | \gamma_j, \tau_j \sim N\left(x_j'\boldsymbol{\beta} + \alpha\gamma_j, \frac{\kappa^2}{\tau_j}\right)$$

$$\gamma_j | \tau_j \sim TN\left(0, \frac{1}{\tau_j}; (0, \infty)\right)$$

$$\tau_j \sim \Gamma\left(\frac{\nu}{2}, \frac{\nu}{2}\right)$$

where, $TN$ denotes the truncated normal distribution, $\alpha = \sigma\delta_\lambda$ and $\kappa^2 = \sigma^2(1 - \delta_\lambda^2)$. Then, regarding $\boldsymbol{\gamma}, \boldsymbol{\tau}$ and $\boldsymbol{z}_j$ are as missing data, the complete data log likelihood function for $(y, \boldsymbol{\gamma}, \boldsymbol{\tau}, \boldsymbol{z}_j)$ given $\boldsymbol{X}$ can be written as

$$\ell_c(\Theta|y, \boldsymbol{\gamma}, \boldsymbol{\tau}, \boldsymbol{z}_j) = \sum_{j=1}^{n}\sum_{i=1}^{g} z_{ij}\left\{\log w_i - \log \pi - \frac{1}{2}\log \kappa_i^2 + \frac{\nu_i}{2}\log\left(\frac{\nu_i}{2}\right) - \log\left(\Gamma\left(\frac{\nu_i}{2}\right)\right) + \frac{\nu_i}{2}\log \tau_j\right.$$

$$-\frac{\nu_i \tau_j}{2} - \frac{(y_j - x_j'\boldsymbol{\beta}_i - \alpha_i \gamma_j)^2}{2\kappa_i^2/\tau_j} - \frac{\gamma_j^2}{2/\tau_j} \right\} \tag{6}$$

where $X = (x_1, \cdots, x_n)'$, $y = (y_1, \cdots, y_n)$, $\gamma = (\gamma_1, \cdots, \gamma_n)$, $\tau = (\tau_1, \cdots, \tau_n)$ and $z_j = (z_{1j}, \cdots, z_{gj})$. Further, based on the theory of the EM algorithm, the conditional expectation of the complete data log-likelihood function given the observed data and the current parameter estimate $\widehat{\Theta}^{(k)}$ should be calculated. That is, we have to find the following conditional expectation

$$E(\ell_c(\Theta|y, \gamma, \tau, z_j)|y_j) = \sum_{j=1}^{n} \sum_{i=1}^{g} E(Z_{ij}|y_j) \left\{ \log w_i - \frac{1}{2}\log \kappa_i^2 + \frac{\nu_i}{2}\log\left(\frac{\nu_i}{2}\right) - \log\left(\Gamma\left(\frac{\nu_i}{2}\right)\right) \right\}$$
$$+ \frac{\nu_i}{2} E(Z_{ij} \log \tau_j | y_j) - \frac{\nu_i}{2} E(Z_{ij}\tau_j|y_j) - \frac{E(Z_{ij}\tau_j|y_j)(y_j - x_j'\boldsymbol{\beta}_i)^2}{2\kappa_i^2}$$
$$+ \frac{\alpha_i E(Z_{ij}\gamma_j\tau_j|y_j)(y_j - x_j'\boldsymbol{\beta}_i)}{\kappa_i^2} - \frac{\alpha_i^2 E(Z_{ij}\gamma_j^2\tau_j|y_j)}{2\kappa_i^2} \tag{7}$$

To get this conditional expectation the following expectations should be obtained: $E(Z_{ij}|y_j, \widehat{\Theta}^{(k)})$, $E(Z_{ij}\tau_j|y_j, \widehat{\Theta}^{(k)})$, $E(Z_{ij}\gamma_j\tau_j|y_j, \widehat{\Theta}^{(k)})$, $E(Z_{ij}\gamma_j^2\tau_j|y_j, \widehat{\Theta}^{(k)})$ and $E(Z_{ij}\log(\tau_j)|y_j, \widehat{\Theta}^{(k)})$. After some straight forward algebra we get the following expressions for these expectations

$$\hat{z}_{ij}^{(k)} = E(Z_{ij}|y_j, \widehat{\Theta}^{(k)}) = \frac{\widehat{w}_i^{(k)} f_i\left(y_j | x_j'\widehat{\boldsymbol{\beta}}_i^{(k)}, \hat{\sigma}_i^{2(k)}, \hat{\lambda}_i^{(k)}, \hat{\nu}_i^{(k)}\right)}{f(y_j|x_j, \widehat{\Theta}^{(k)})}, \tag{8}$$

$$\hat{s}_{1ij}^{(k)} = E(Z_{ij}\tau_j|y_j, \widehat{\Theta}^{(k)}) = \hat{z}_{ij}^{(k)} \left(\frac{\hat{\nu}_i^{(k)} + 1}{\hat{\eta}_{ij}^{2(k)} + \hat{\nu}_i^{(k)}}\right) \frac{T_{\hat{\nu}_i^{(k)}+3}\left(\widehat{M}_{ij}^{(k)} \sqrt{\frac{\hat{\nu}_i^{(k)}+3}{\hat{\nu}_i^{(k)}+1}}\right)}{T_{\hat{\nu}_i^{(k)}+1}\left(\widehat{M}_{ij}^{(k)}\right)}, \tag{9}$$

$$\hat{s}_{2ij}^{(k)} = E(Z_{ij}\gamma_j\tau_j|y_j, \widehat{\Theta}^{(k)}) = \frac{1}{\hat{\sigma}_i^{(k)}} \hat{\delta}_{\lambda_i}^{(k)} \left(y_j - x_j'\widehat{\boldsymbol{\beta}}_i^{(k)}\right) \hat{s}_{1ij}^{(k)}$$
$$+ \hat{z}_{ij}^{(k)} \frac{\sqrt{1 - \hat{\delta}_{\lambda_i}^{2(k)}}}{\pi \hat{\sigma}_i^{(k)} \hat{f}(y_j)^{(k)}} \left(\frac{\hat{\eta}_{ij}^{2(k)}}{\hat{\nu}_i^{(k)}\left(1 - \hat{\delta}_{\lambda_i}^{2(k)}\right)} + 1\right)^{-\left(\frac{\hat{\nu}_i^{(k)}}{2} + 1\right)}, \tag{10}$$

$$\hat{s}_{3ij}^{(k)} = E(Z_{ij}\gamma_j^2\tau_j|y_j, \widehat{\Theta}^{(k)})$$
$$= \hat{\delta}_{\lambda_i}^{2(k)} \left(\frac{y_j - x_j'\widehat{\boldsymbol{\beta}}_i^{(k)}}{\hat{\sigma}_i^{(k)}}\right)^2 \hat{s}_{1ij}^{(k)} + \hat{z}_{ij}^{(k)} \left\{ \left(1 - \hat{\delta}_{\lambda_i}^{2(k)}\right) \right.$$
$$\left. + \frac{\hat{\delta}_{\lambda_i}^{(k)} \left(y_j - x_j'\widehat{\boldsymbol{\beta}}_i^{(k)}\right) \sqrt{1 - \hat{\delta}_{\lambda_i}^{2(k)}}}{\pi \hat{\sigma}_i^{2(k)} \hat{f}(y_j)^{(k)}} \left(\frac{\hat{\eta}_{ij}^{2(k)}}{\hat{\nu}_i^{(k)}\left(1 - \hat{\delta}_{\lambda_i}^{2(k)}\right)} + 1\right)^{-\left(\frac{\hat{\nu}_i^{(k)}}{2} + 1\right)} \right\} \tag{11}$$

$$\hat{s}_{4ij}^{(k)} = E(Z_{ij}\log(\tau_j)|y_j, \widehat{\Theta}^{(k)})$$
$$= \hat{z}_{ij}^{(k)} \left\{ DG\left(\frac{\hat{\nu}_i^{(k)} + 1}{2}\right) - \log\left(\frac{\hat{\eta}_{ij}^{2(k)} + \hat{\nu}_i^{(k)}}{2}\right) \right.$$

$$+\left(\frac{\hat{v}_i^{(k)}+1}{\hat{\eta}_{ij}^{2(k)}+\hat{v}_i^{(k)}}\right)\left(\frac{T_{\hat{v}_i^{(k)}+3}\left(\hat{\lambda}_i^{(k)}\hat{\eta}_{ij}^{(k)}\sqrt{\frac{\hat{v}_i^{(k)}+3}{\hat{\eta}_{ij}^{2(k)}+\hat{v}_i^{(k)}}}\right)}{T_{\hat{v}_i^{(k)}+1}\left(\hat{\lambda}_i^{(k)}\hat{\eta}_{ij}^{(k)}\sqrt{\frac{\hat{v}_i^{(k)}+1}{\hat{\eta}_{ij}^{2(k)}+\hat{v}_i^{(k)}}}\right)}-1\right)$$

$$+\frac{\hat{\lambda}_i^{(k)}\hat{\eta}_{ij}^{(k)}\left(\hat{\eta}_{ij}^{2(k)}-1\right)}{\sqrt{\left(\hat{v}_i^{(k)}+1\right)\left(\hat{\eta}_{ij}^{2(k)}+\hat{v}_i^{(k)}\right)^3}}\frac{t_{\hat{v}_i^{(k)}+1}\left(\hat{\lambda}_i^{(k)}\hat{\eta}_{ij}^{(k)}\sqrt{\frac{\hat{v}_i^{(k)}+1}{\hat{\eta}_{ij}^{2(k)}+\hat{v}_i^{(k)}}}\right)}{T_{\hat{v}_i^{(k)}+1}\left(\hat{\lambda}_i^{(k)}\hat{\eta}_{ij}^{(k)}\sqrt{\frac{\hat{v}_i^{(k)}+1}{\hat{\eta}_{ij}^{2(k)}+\hat{v}_i^{(k)}}}\right)}$$

$$\left.+\frac{1}{T_{\hat{v}_i^{(k)}+1}\left(\hat{\lambda}_i^{(k)}\hat{\eta}_{ij}^{(k)}\sqrt{\frac{\hat{v}_i^{(k)}+1}{\hat{\eta}_{ij}^{2(k)}+\hat{v}_i^{(k)}}}\right)}\int_{-\infty}^{\widehat{M}_{ij}^{(k)}}g_{\hat{v}_i^{(k)}}(x)t_{\hat{v}_i^{(k)}+1}(x)dx\right\} \quad (12)$$

where

$$\hat{\eta}_{ij}^{(k)}=\frac{\left(y_j-x_j'\hat{\boldsymbol{\beta}}_i^{(k)}\right)}{\hat{\sigma}_i^{(k)}},\delta_{\hat{\lambda}_i}^{(k)}=\frac{\hat{\lambda}_i^{(k)}}{\sqrt{1+\hat{\lambda}_i^{2(k)}}},$$

$$\widehat{M}_{ij}^{(k)}=\hat{\lambda}_i^{(k)}\hat{\eta}_{ij}^{(k)}\sqrt{\frac{\hat{v}_i^{(k)}}{\hat{\eta}_{ij}^{2(k)}+\hat{v}_i^{(k)}}},\hat{f}(y_j)^{(k)}=\sum_{i=1}^g\hat{w}_i^{(k)}\frac{2}{\hat{\sigma}_i^{(k)}}t_{\hat{v}_i^{(k)}}\left(\hat{\eta}_{ij}^{(k)}\right)T_{\hat{v}_i^{(k)}+1}\left(\widehat{M}_{ij}^{(k)}\right).$$

Then, the EM algorithm to obtain the parameter estimates for the mixture regression model based on the skew *t* distribution can be given as follows.

**EM algorithm:**

**1.** Take initial parameter estimates $\boldsymbol{\Theta}^{(0)}$.

**2.** E step: To proceed an E step, we have to find the conditional expectation of the complete data log likelihood function given the current parameter values $\boldsymbol{\Theta}^{(k)}$. This can be done by computing the conditional expectations $\hat{z}_{ij}^{(k)}, \hat{s}_{1ij}^{(k)}, \hat{s}_{2ij}^{(k)}, \hat{s}_{3ij}^{(k)}$ and $\hat{s}_{4ij}^{(k)}$ for $k=0,1,2,...$. After finding these conditional expectations we get the following objective function to be maximized at M step of the EM algorithm

$$Q(\boldsymbol{\Theta}|\widehat{\boldsymbol{\Theta}}^{(k)})=\sum_{j=1}^n\sum_{i=1}^g\hat{z}_{ij}^{(k)}\left\{\log w_i-\frac{1}{2}\log\kappa_i^2+\frac{v_i}{2}\log\left(\frac{v_i}{2}\right)-\log\left(\Gamma\left(\frac{v_i}{2}\right)\right)\right\}$$

$$+\frac{v_i}{2}\hat{s}_{4ij}^{(k)}-\frac{v_i}{2}\hat{s}_{1ij}^{(k)}-\frac{\hat{s}_{1ij}^{(k)}\left(y_j-x_j'\boldsymbol{\beta}_i\right)^2}{2\kappa_i^2}+\frac{\alpha_i\hat{s}_{2ij}^{(k)}(y_j-x_j'\boldsymbol{\beta}_i)}{\kappa_i^2}-\frac{\alpha_i^2\hat{s}_{3ij}^{(k)}}{2\kappa_i^2}. \quad (13)$$

**3. M step 1:** Maximize the $Q(\Theta|\widehat{\Theta}^{(k)})$ with respect to the unknown parameters $(w_i, \boldsymbol{\beta}_i, \sigma_i^2)$, assuming that $(\lambda_i, \nu_i)$ are fixed, to obtain $(k+1)$th values for the parameter $(w_i, \boldsymbol{\beta}_i, \sigma_i^2)$. This maximization gives

$$w_i^{(k+1)} = \frac{\sum_{j=1}^n \hat{z}_{ij}^{(k)}}{n}, \tag{14}$$

$$\widehat{\boldsymbol{\beta}}_i^{(k+1)} = \left(\sum_{j=1}^n x_j x_j' \hat{s}_{1ij}^{(k)}\right)^{-1} \left(\sum_{j=1}^n \left(y_j \hat{s}_{1ij}^{(k)} - \hat{\alpha}_i^{(k)} \hat{s}_{2ij}^{(k)}\right) x_j\right), \tag{15}$$

$$\hat{\alpha}_i^{(k+1)} = \frac{\sum_{j=1}^n \hat{s}_{2ij}^{(k)} \left(y_j - x_j'\widehat{\boldsymbol{\beta}}_i^{(k)}\right)}{\sum_{j=1}^n \hat{s}_{3ij}^{(k)}}, \tag{16}$$

$$\hat{\kappa}_i^{(k+1)} = \frac{\sum_{j=1}^n \left(\hat{s}_{1ij}^{(k)} \left(y_j - x_j'\widehat{\boldsymbol{\beta}}_i^{(k)}\right)^2 - 2\hat{\alpha}_i^{(k)} \hat{s}_{2ij}^{(k)} \left(y_j - x_j'\widehat{\boldsymbol{\beta}}_i^{(k)}\right) + \hat{\alpha}_i^{2(k)} \hat{s}_{3ij}^{(k)}\right)}{\sum_{j=1}^n \hat{z}_{ij}^{(k)}}, \tag{17}$$

$$\hat{\sigma}_i^{(k+1)} = \hat{\kappa}_i^{(k+1)} + \hat{\alpha}_i^{(k+1)}. \tag{18}$$

**4. M step 2:** Using the new values for $(w_i, \boldsymbol{\beta}_i, \sigma_i^2)$ gained in M step 1 solve the following equations to obtain new estimates for the parameters $(\lambda_i, \nu_i)$

$$\delta_{\lambda_i}(1-\delta_{\lambda_i}^2)\sum_{j=1}^n \hat{z}_{ij}^{(k)} - \delta_{\lambda_i}\left(\sum_{j=1}^n \hat{s}_{1ij}^{(k)} \frac{\left(y_j - x_j'\boldsymbol{\beta}_i^{(k+1)}\right)^2}{\sigma_i^{2(k+1)}} + \sum_{j=1}^n \hat{s}_{3ij}^{(k)}\right)$$
$$+ (1+\delta_{\lambda_i}^2)\sum_{j=1}^n \hat{s}_{2ij}^{(k)} \frac{\left(y_j - x_j'\boldsymbol{\beta}_i^{(k+1)}\right)}{\sigma_i^{(k+1)}} = 0, \tag{19}$$

$$\log\left(\frac{\nu_i}{2}\right) + 1 - DG\left(\frac{\nu_i}{2}\right) + \frac{\frac{1}{n}\sum_{j=1}^n \left(\hat{s}_{4j}^{(k)} - \hat{s}_{1j}^{(k)}\right)}{\sum_{j=1}^n \hat{z}_{ij}^{(k)}} = 0. \tag{20}$$

**5.** Repeat E and M steps until the convergence criteria $\|\Theta^{(k+1)} - \Theta^{(k)}\| < \epsilon$ is satisfied.

Note that to simplify the computation of $\hat{\lambda}_i^{(k+1)}$ we will use the following estimate in the simulation study and real data example

$$\hat{\lambda}_i^{(k+1)} = \hat{\delta}_{\lambda_i}^{(k+1)} / \sqrt{1 - \hat{\delta}_{\lambda_i}^{2(k+1)}},$$

where $\hat{\delta}_{\lambda_i}^{(k+1)} = \hat{\alpha}_i^{(k+1)}/\hat{\sigma}_i^{(k+1)}$.

## 4. Simulation Study

In this section, we will give a simulation study to show the performance of the proposed estimator obtained from skew *t* (MixregST) and we also compare the other estimators obtained from normal (MixregN), *t* (Mixregt) and skew normal (MixregSN) distributions in terms of bias and mean square error (MSE).

We generate the data $\{(x_{1i}, x_{2i}, y_i), i = 1, \ldots, n\}$ from the following two component mixture regression models (Bai et al. (2012))

$$Y = \begin{cases} 0 + X_1 + X_2 + \epsilon_1, & Z = 1, \\ 0 - X_1 - X_2 + \epsilon_2, & Z = 2, \end{cases}$$

where $P(Z = 1) = 0.25 = w_1, X_1 \sim N(0,1), X_2 \sim N(0,1)$, the errors are i.i.d. Furthermore, the model coefficients are $\beta_1 = (\beta_{10}, \beta_{11}, \beta_{12})' = (0,1,1)'$ and $\beta_2 = (\beta_{20}, \beta_{21}, \beta_{22})' = (0,-1,-1)'$.

We take the following error distributions:

Case I: $\epsilon \sim N(0,1)$, standard normal distribution.
Case II: $\epsilon \sim t_3(0,1)$, $t$ distribution with the degrees of freedom 3.
Case III: $\epsilon \sim 0.95N(0,1) + 0.05N(0,25)$, contaminated normal distribution.
Case IV: $\epsilon \sim ST(0,1,0.5,3)$, skew $t$ distribution.
Case V: $\epsilon \sim N(0,1)$, standard normal distribution with %5 outliers, $X_1 = 20, X_2 = 20$ and $Y = 100$.

We use the Case I to compare the estimators with the traditional MLE (MixregN) when the error terms have the normal distribution and there are no outliers. Case II is the example for the heavy-tailed error distribution. The distribution given in Case III is to create outliers. This distribution is often considered in literature as an outlier model. Case IV is to examine the behavior of the estimators when the error term is skewed and heavy-tailed. Case V is considered to test the performances of the estimators to deal with the high leverage points. In this case %5 of the observations are replaced by $X_1 = 20, X_2 = 20$ and $Y = 100$. In the simulation study, the sample sizes are taken as 200 and 400 and the number of replicates is 500. The simulation study and real data example are conducted using MATLAB 2013a.

Table 1 and 2 show the MSEs and the biases of the parameter estimates. We can observe from the results of the simulation study that the MixregN has the best result in Case I. On the other hand, the other estimators Mixregt, MixregSN and MixregST have similar performances when the errors have normal distribution. In Case II, Mixregt performs best, as expected. Also, MixregST has the lower bias and MSE values than the MixregN and MixregSN for almost all the cases. For the Case III, MixregN and MixregSN are drastically affected by the contamination. However, Mixregt and MixregST perform better than the other estimators and Mixregt is comparable with the MixregST. Similarly, MixregN and MixregSN have the worst performance and Mixregt and MixregST have similar performance in Case IV. Finally, in the outlier case all estimators are affected by the outliers. However, Mixregt and MixregST have the lowest bias and MSE values in almost all cases. In summary, concerning all the estimators the Mixregt and MixregST are resistant to the skewness and the heavy tailedness in the data, and they behave better than MixregN and MixregSN in case of outliers in $X$ direction.

**Table 1**
MSE (bias) values of estimates for $n = 200$.

| | MixregN | Mixregt | MixregSN | MixregST |
|---|---|---|---|---|
| | Case I: $\epsilon \sim N(0,1)$ | | | |
| $\beta_{10}: 0$ | 0.0456 (0.0150) | 0.0587 (0.0134) | 0.1726 (-0.3560) | 0.1317 (0.2306) |
| $\beta_{20}: 0$ | 0.0090 (0.0019) | 0.0098 (0.0039) | 0.1447 (-0.3678) | 0.0575 (-0.2084) |
| $\beta_{11}: 1$ | 0.0348 (-0.0013) | 0.0495 (-0.0064) | 0.0349 (-0.0016) | 0.0546 (-0.0036) |
| $\beta_{21}: -1$ | 0.0085 (-0.0004) | 0.0103 (0.0031) | 0.0085 (-0.0004) | 0.0118 (0.0212) |
| $\beta_{12}: 1$ | 0.0401 (-0.0243) | 0.0483 (-0.0308) | 0.0401 (-0.0242) | 0.0617 (-0.0296) |
| $\beta_{22}: -1$ | 0.0089 (-0.0004) | 0.0107 (0.0024) | 0.0089 (-0.0062) | 0.0125 (0.0201) |
| $w_1: 0.25$ | 0.0021 (0.0079) | 0.0023 (0.0059) | 0.0021 (0.0079) | 0.0035 (-0.0063) |
| | Case II: $\epsilon \sim t_3$ | | | |
| $\beta_{10}: 0$ | 11.5674 (-0.2939) | 0.0930 (-0.0121) | 11.6586 (-0.9305) | 0.3151 (0.2406) |
| $\beta_{20}: 0$ | 1.2217 (0.0796) | 0.0136 (-0.0050) | 1.3914 (-0.5527) | 0.1397 (-0.3327) |
| $\beta_{11}: 1$ | 7.6108 (0.4273) | 0.0959 (-0.0180) | 7.6526 (0.3704) | 0.1415 (0.0036) |
| $\beta_{21}: -1$ | 1.2984 (-0.0331) | 0.0145 (-0.0064) | 1.2011 (0.0192) | 0.0171 (0.0259) |
| $\beta_{12}: 1$ | 8.2789 (0.1660) | 0.0981 (0.0027) | 8.2956 (0.2624) | 0.1678 (0.0282) |
| $\beta_{22}: -1$ | 1.9409 (-0.0331) | 0.0137 (-0.0031) | 1.6075 (0.0762) | 0.0167 (0.0283) |
| $w_1: 0.25$ | 0.0226 (-0.0372) | 0.0033 (0.0112) | 0.0214 (-0.0352) | 0.0055 (-0.0067) |
| | Case III: $\epsilon \sim 0.95N(0,1) + 0.05N(0,25)$ | | | |
| $\beta_{10}: 0$ | 6.0158 (-0.0052) | 0.0634 (-0.0062) | 6.1249 (-0.6206) | 0.1517 (0.2053) |
| $\beta_{20}: 0$ | 0.6299 (0.0054) | 0.0118 (-0.0080) | 0.6282 (-0.5670) | 0.0911 (-0.2711) |
| $\beta_{11}: 1$ | 4.5781 (0.2371) | 0.0599 (0.0078) | 4.8849 (0.2067) | 0.0727 (0.0119) |
| $\beta_{21}: -1$ | 0.2236 (0.0418) | 0.0106 (-0.0068) | 0.1302 (0.0649) | 0.0124 (0.0155) |
| $\beta_{12}: 1$ | 2.9126 (-0.0271) | 0.0620 (0.0021) | 2.7706 (0.0830) | 0.0774 (0.0192) |
| $\beta_{22}: -1$ | 0.1607 (0.0418) | 0.0090 (0.0033) | 0.0614 (0.0778) | 0.0108 (0.0250) |
| $w_1: 0.25$ | 0.0167 (-0.0472) | 0.0026 (0.0039) | 0.0136 (-0.0526) | 0.0034 (-0.0098) |
| | Case IV: $\epsilon \sim ST(0,1,0.5,3)$ | | | |
| $\beta_{10}: 0$ | 8.4499 (1.0601) | 0.2783 (0.4422) | 6.1264 (0.3167) | 0.9691 (0.7550) |
| $\beta_{20}: 0$ | 0.3472 (0.4787) | 0.1524 (0.3759) | 0.1323 (-0.0886) | 0.0231 (0.0590) |
| $\beta_{11}: 1$ | 2.9291 (0.2448) | 0.0851 (-0.0296) | 2.7053 (0.2225) | 0.1605 (-0.0107) |
| $\beta_{21}: -1$ | 0.0600 (0.0432) | 0.0120 (-0.0133) | 0.0540 (0.0381) | 0.0146 (0.0230) |
| $\beta_{12}: 1$ | 5.9774 (-0.1412) | 0.0862 (-0.0195) | 5.6460 (-0.0863) | 0.1911 (0.0005) |
| $\beta_{22}: -1$ | 0.0789 (0.0432) | 0.0115 (-0.0029) | 0.0731 (0.0715) | 0.0154 (0.0336) |
| $w_1: 0.25$ | 0.0125 (-0.0296) | 0.0033 (0.0118) | 0.0116 (-0.0260) | 0.0050 (-0.0156) |
| | Case V: $\epsilon \sim N(0,1)$ (%5 outliers) | | | |
| $\beta_{10}: 0$ | 2.2247 (0.1553) | 1.3245 (0.1820) | 2.5926 (-0.4879) | 5.9114 (2.1745) |
| $\beta_{20}: 0$ | 0.0146 (0.0111) | 0.0106 (0.0072) | 0.2401 (-0.4728) | 0.0392 (-0.1678) |
| $\beta_{11}: 1$ | 3.2773 (1.5211) | 2.8341 (1.5030) | 3.3162 (1.5107) | 2.6095 (1.4250) |
| $\beta_{21}: -1$ | 0.0833 (0.2528) | 0.0234 (0.1077) | 0.0826 (0.2519) | 0.0296 (0.1283) |
| $\beta_{12}: 1$ | 3.1162 (1.4674) | 2.7897 (1.4869) | 3.2436 (1.4870) | 2.7237 (1.4655) |
| $\beta_{22}: -1$ | 0.0798 (0.2528) | 0.0225 (0.1055) | 0.0786 (0.2472) | 0.0281 (0.1244) |
| $w_1: 0.25$ | 0.0093 (-0.0937) | 0.0061 (-0.0751) | 0.0094 (-0.0939) | 0.0112 (-0.1029) |

**Table 2**
MSE (bias) values of estimates for $n = 400$.

| | MixregN | Mixregt | MixregSN | MixregST |
|---|---|---|---|---|
| | Case I: $\epsilon \sim N(0,1)$ | | | |
| $\beta_{10}: 0$ | 0.0203 (0.0088) | 0.0265 (0.0114) | 0.1564 (-0.3687) | 0.0782 (0.2081) |
| $\beta_{20}: 0$ | 0.0043 (0.0044) | 0.0050 (0.0058) | 0.1427 (-0.3716) | 0.0565 (-0.2211) |
| $\beta_{11}: 1$ | 0.0149 (0.0008) | 0.0192 (-0.0019) | 0.0149 (0.0009) | 0.0227 (0.0034) |
| $\beta_{21}: -1$ | 0.0040 (-0.0028) | 0.0048 (0.0016) | 0.0040 (-0.0027) | 0.0057 (0.0197) |
| $\beta_{12}: 1$ | 0.0160 (-0.0100) | 0.0213 (-0.0185) | 0.0161 (-0.0100) | 0.0245 (-0.0110) |
| $\beta_{22}: -1$ | 0.0044 (-0.0028) | 0.0053 (0.0070) | 0.0044 (0.0010) | 0.0065 (0.0247) |
| $w_1: 0.25$ | 0.0012 (0.0035) | 0.0013 (0.0006) | 0.0012 (0.0035) | 0.0018 (-0.0123) |
| | Case II: $\epsilon \sim t_3$ | | | |
| $\beta_{10}: 0$ | 14.3296 (-0.3312) | 0.0365 (-0.0137) | 14.2254 (-0.9701) | 0.0830 (0.1669) |
| $\beta_{20}: 0$ | 0.6052 (0.0125) | 0.0066 (-0.0066) | 0.6330 (-0.6752) | 0.1411 (-0.3601) |
| $\beta_{11}: 1$ | 10.6597 (0.4839) | 0.0321 (-0.0054) | 10.1135 (0.4010) | 0.0427 (0.0239) |
| $\beta_{21}: -1$ | 0.5987 (0.0527) | 0.0068 (-0.0052) | 0.1809 (0.0921) | 0.0083 (0.0272) |
| $\beta_{12}: 1$ | 12.1779 (0.3384) | 0.0334 (0.0052) | 11.8293 (0.5888) | 0.0421 (0.0288) |
| $\beta_{22}: -1$ | 1.5732 (0.0527) | 0.0062 (-0.0041) | 0.9058 (0.0454) | 0.0078 (0.0273) |
| $w_1: 0.25$ | 0.0161 (-0.0602) | 0.0014 (0.0049) | 0.0143 (-0.0591) | 0.0020 (-0.0134) |
| | Case III: $\epsilon \sim 0.95N(0,1) + 0.05N(0,25)$ | | | |
| $\beta_{10}: 0$ | 4.6683 (-0.0431) | 0.0287 (0.0004) | 5.1651 (-0.7278) | 0.0729 (0.1830) |
| $\beta_{20}: 0$ | 0.0088 (0.0037) | 0.0056 (-0.0012) | 0.3555 (-0.5817) | 0.0848 (-0.2769) |
| $\beta_{11}: 1$ | 4.2093 (0.1003) | 0.0229 (0.0038) | 4.2202 (0.0962) | 0.0278 (0.0214) |
| $\beta_{21}: -1$ | 0.0313 (0.0872) | 0.0053 (0.0024) | 0.0319 (0.0875) | 0.0066 (0.0243) |
| $\beta_{12}: 1$ | 3.2445 (0.1817) | 0.0251 (0.0166) | 3.1090 (0.1611) | 0.0327 (0.0303) |
| $\beta_{22}: -1$ | 0.0328 (0.0872) | 0.0054 (0.0064) | 0.0325 (0.0878) | 0.0069 (0.0292) |
| $w_1: 0.25$ | 0.0093 (-0.0572) | 0.0014 (-0.0020) | 0.0093 (-0.0570) | 0.0019 (-0.0160) |
| | Case IV: $\epsilon \sim ST(0,1,0.5,3)$ | | | |
| $\beta_{10}: 0$ | 7.8868 (0.9770) | 0.2082 (0.4344) | 5.1906 (0.0754) | 0.4395 (0.6371) |
| $\beta_{20}: 0$ | 0.2110 (0.4604) | 0.1461 (0.3853) | 0.0373 (-0.1476) | 0.0105 (0.0455) |
| $\beta_{11}: 1$ | 5.0109 (0.1370) | 0.0247 (-0.0192) | 5.6461 (0.1695) | 0.0400 (0.0140) |
| $\beta_{21}: -1$ | 0.0280 (0.0717) | 0.0053 (-0.0065) | 0.0259 (0.0686) | 0.0066 (0.0263) |
| $\beta_{12}: 1$ | 6.6126 (0.3814) | 0.0301 (-0.0120) | 7.0245 (0.3604) | 0.0485 (0.0140) |
| $\beta_{22}: -1$ | 0.0308 (0.0717) | 0.0049 (-0.0044) | 0.0276 (0.0691) | 0.0069 (0.0290) |
| $w_1: 0.25$ | 0.0081 (-0.0459) | 0.0014 (0.0040) | 0.0073 (-0.0436) | 0.0021 (-0.0168) |
| | Case V: $\epsilon \sim N(0,1)$ (%5 outliers) | | | |
| $\beta_{10}: 0$ | 1.5208 (0.2485) | 1.0975 (0.2305) | 1.6056 (-0.3105) | 6.9413 (2.5194) |
| $\beta_{20}: 0$ | 0.0094 (0.0158) | 0.0059 (0.0056) | 0.2483 (-0.4883) | 0.0419 (-0.1880) |
| $\beta_{11}: 1$ | 2.6872 (1.4449) | 2.4663 (1.4533) | 2.6444 (1.4307) | 2.3970 (1.4530) |
| $\beta_{21}: -1$ | 0.0783 (0.2591) | 0.0175 (0.1066) | 0.0770 (0.2572) | 0.0239 (0.1284) |
| $\beta_{12}: 1$ | 2.9720 (1.5383) | 2.7078 (1.5341) | 3.0209 (1.5560) | 2.3044 (1.4204) |
| $\beta_{22}: -1$ | 0.0813 (0.2591) | 0.0176 (0.1072) | 0.0810 (0.2639) | 0.0230 (0.1279) |
| $w_1: 0.25$ | 0.0098 (-0.0974) | 0.0069 (-0.0814) | 0.0098 (-0.0976) | 0.0138 (-0.1159) |

## 5. Real Data Example

In this section we analyze the tone perception data set (Cohen (1984)) to further illustrate the performance of the mixture regression estimates based on the skew *t* distribution on a real data set. In the tone perception experiment of Cohen (1984), a pure fundamental tone was played to a trained musician. Also, electronically obtained overtones were added which were determined by a stretching ratio. This ratio is between the adjusted tone and the fundamental tone. In the experiment, 150 trials were performed by the same musicians. The aim of this experiment was to find out how the tuning ratio affects the perception of the tone and to decide if either of two musical perception theories was reasonable (see Cohen (1984) for more detail). This data set has also been analyzed by Yao et al. (2014) and Song et al. (2014) to test the performance of the mixture regression estimates based on the *t* and Laplace distributions, respectively. Figure 1 shows the scatter plot and the histogram of the perceived tone ratio. From these plots it is clear that there are two groups in the data and it also shows the non-normality.

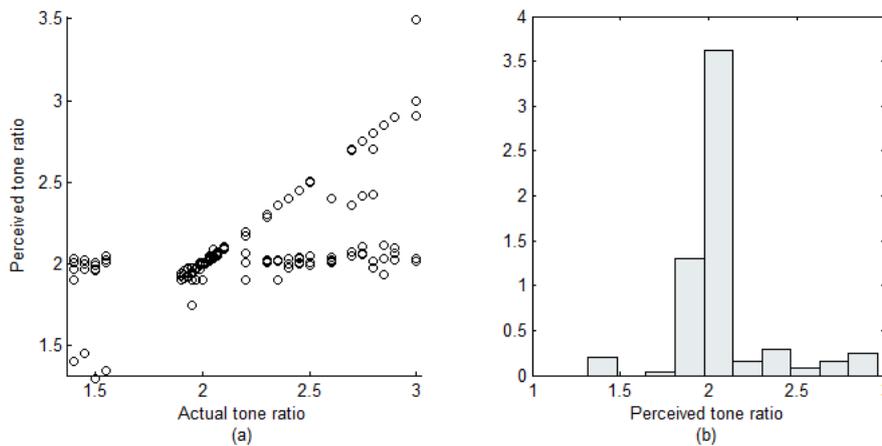

Figure 1. (a) The scatter plot of the data. (b) Histogram of the perceived tone ratio.

We use this data set to compare the performances of the estimators in the case of with and without outliers. We present the scatter plots with the fitted regression lines obtained from MixregN, Mixregt, MixregSN and MixregST procedures in Figure 2 for the tone perception data set. Also, we summary the ML estimates and some information criterions in Table 3. Note that in real data example we assume that in both groups the degrees of freedom equals to 2. We try other values of degrees of freedom and get the similar results. We observe that MixregST has the best fit than the other mixture regression models in terms of the Akaike information criterion (AIC) (Akaike (1973)) and the Bayesian information criterion (BIC) (Schwarz (1978)) values.

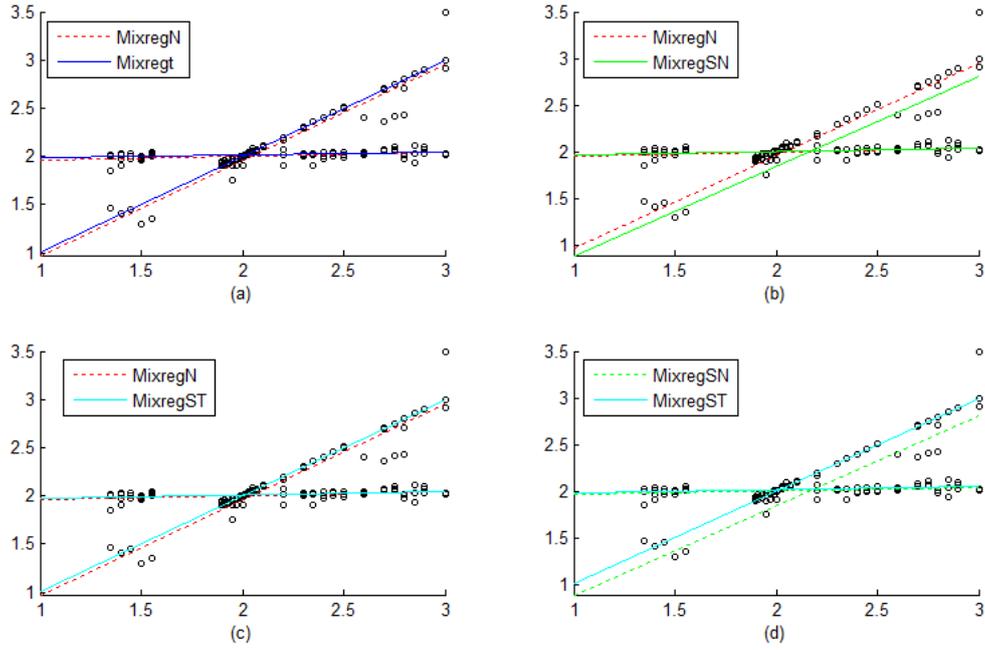

Figure 2. Fitted mixture regression lines for the tone perception data set. (a): dashed line- MixregN, solid line-Mixregt, (b): dashed line- MixregN, solid line-MixregSN, (c): dashed line- MixregN, solid line-MixregST, (d): dashed line-MixregSN, solid line-MixregST.

**Table 3.** ML estimates and some information criterions for fitting different mixture regression models to the tone perception data set

|  | MixregN | Mixregt | MixregSN | MixregST |
|---|---|---|---|---|
| $\hat{\beta}_{10}$ | 1.9164 | 1.9586 | 1.9171 | 1.9491 |
| $\hat{\beta}_{20}$ | 0.0425 | 0.0264 | 0.0424 | 0.0318 |
| $\hat{\beta}_{11}$ | -0.0193 | 0.0178 | -0.0717 | 0.0054 |
| $\hat{\beta}_{21}$ | 0.9923 | 0.9918 | 0.9604 | 0.9982 |
| $\hat{\sigma}_1$ | 0.0462 | 0.0281 | 0.0463 | 0.0393 |
| $\hat{\sigma}_2$ | 0.1328 | 0.0210 | 0.1883 | 0.0033 |
| $\hat{\lambda}_1$ | - | - | -0.0100 | -0.1666 |
| $\hat{\lambda}_2$ | - | - | 1.7534 | 0.4465 |
| $\hat{w}_1$ | 0.6977 | 0.5518 | 0.7006 | 0.6410 |
| $\ell(\hat{\Theta})$ | 141.1984 | 190.8177 | 140.5585 | **211.7766** |
| AIC | -286.3968 | -367.6354 | -263.1171 | **-405.5532** |
| BIC | -247.3224 | -346.5610 | -236.0213 | **-378.4574** |

Next we add ten pairs of outliers at (0,5). These outliers can be considered as high leverage points. By adding these points we would like to see the performance of the estimators against to the high leverage points. Figure 3 displays the scatter plots of the data set with the fitted regression lines obtained from MixregN, Mixregt, MixregSN and MixregST procedures. We give the ML estimation results in Table 4. We see that MixregN and MixregSN are drastically affected by the high leverage points. On the other hands, the estimators based on the *t* and the skew *t* distributions (Mixregt and MixregST) give fits to the majority of the data without influencing from the high leverage points. Also, MixregST gives best results in terms of information criterion. Note that the estimates including the estimates for skewness parameters with and without outliers are very similar (see Tables 3 and 4)

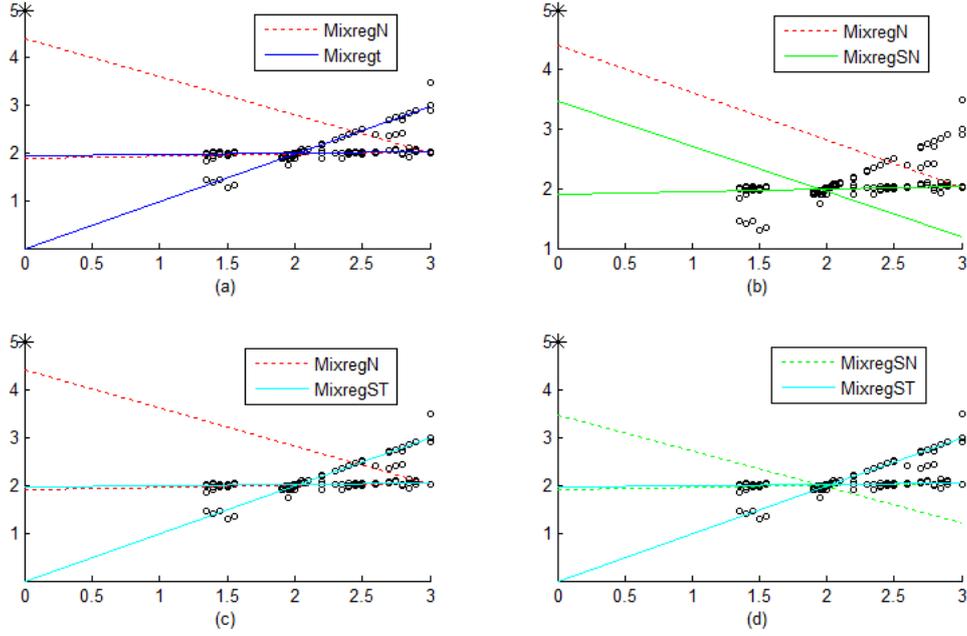

Figure 3. Fitted mixture regression lines with ten outliers at (0,5). (a): dashed line- MixregN, solid line-Mixregt, (b): dashed line- MixregN, solid line-MixregSN, (c): dashed line- MixregN, solid line-MixregST, (d): dashed line-MixregSN, solid line-MixregST.

**Table 4.** ML estimates and some information criterions for fitting different mixture regression models to the tone perception data set with ten outliers at (0,5)

|  | MixregN | Mixregt | MixregSN | MixregST |
|---|---|---|---|---|
| $\hat{\beta}_{10}$ | 1.9058 | 1.9529 | 1.8948 | 1.9553 |
| $\hat{\beta}_{20}$ | 0.0471 | 0.0288 | 0.0478 | 0.0313 |
| $\hat{\beta}_{11}$ | 4.4010 | 0.0251 | 3.4734 | 0.0057 |
| $\hat{\beta}_{21}$ | -0.7954 | 0.9881 | -0.7579 | 0.9981 |
| $\hat{\sigma}_1$ | 0.0506 | 0.0400 | 0.0612 | 0.0542 |
| $\hat{\sigma}_2$ | 0.8591 | 0.0280 | 1.2593 | 0.0031 |
| $\hat{\lambda}_1$ | - | - | -0.2667 | -0.2030 |
| $\hat{\lambda}_2$ | - | - | 1.6770 | 0.4493 |
| $\hat{w}_1$ | 0.7367 | 0.6083 | 0.7382 | 0.6759 |
| $\ell(\hat{\Theta})$ | 54.0997 | 77.5769 | 40.7933 | **109.3612** |
| AIC | -94.1994 | -141.1537 | -63.5867 | **-200.7225** |
| BIC | -72.6732 | -119.6275 | -35.9101 | **-173.0459** |

## 6. Conclusions

In this paper we have proposed a robust mixture regression procedure based on the skew *t* distribution. We have given an EM algorithm to compute the proposed estimators for the mixture regression model. We have given a simulation study to explore the performance of the estimators based on the skew *t* distribution over the estimators obtained from the normal, the *t* and the skew normal distributions. The simulation results confirm that when heavy-tailedness and skewness are present the proposed estimators behave better than the counterparts. We have also given a real data example to further illustrate the capabilty of the proposed estimators dealing with the outliers and/or high leverage points in the data. Likewise, for the real data our proposed estimators show superiorty over the estimators based on normal, *t* and skew normal.

**Appendix**

If a random variable $Y$ has the skew $t$ distribution ($ST(\xi, \sigma^2, \lambda, \nu)$) with the location parameter $\xi \in \mathbb{R}$, scale parameter $\sigma^2 \in (0, \infty)$, skewness parameter $\lambda \in \mathbb{R}$ and degrees of freedom $\nu$, it has the following stochastic representation (Azzalini and Capitaino (2003))

$$Y = \xi + \sigma \frac{Z}{\sqrt{\tau}}, Z \sim SN(\lambda), \tau \sim \Gamma(\nu/2, \nu/2),$$

where $Z$ and $\tau$ are independent and $SN$ shows the skew normal distribution, respectively. Also we can further give the following stochastic representation for $Z$, which has already given by Azzalini (1986, p.201) and Henze (1986, Theorem 1)

$$Z = \delta_\lambda |U_1| + \sqrt{1 - \delta_\lambda^2} U_2,$$

where $U_1$ and $U_2$ are independent standard normal random variables and $|U_1|$ will have truncated normal distribution. This stochastic representation can be used to get the following conditional distributions

$$Y|\gamma, \tau \sim N\left(\xi + \sigma \delta_\lambda \gamma, \frac{1 - \delta_\lambda^2}{\tau} \sigma^2\right),$$

$$\gamma|\tau \sim TN\left(0, \frac{1}{\tau}; (0, \infty)\right), \quad \tau \sim \Gamma(\nu/2, \nu/2).$$

These conditional distributions will help us to conduct the steps of the EM algorithm. By Proposition 2 of Lin et al. (2007) we can have the following conditional expectations for $\tau$, $\gamma\tau$, $\gamma^2\tau$ and $\log(\tau)$ given $Y = y$

$$\boldsymbol{i}) \; E(\tau|y) = \left(\frac{\nu + 1}{\eta^2 + \nu}\right) \frac{T_{\nu+3}\left(M \sqrt{\frac{\nu+3}{\nu+1}}\right)}{T_{\nu+1}(M)},$$

$$\boldsymbol{ii}) \; E(\gamma\tau|y) = \delta_\lambda \frac{(y - \xi)}{\sigma} E(\tau|y) + \frac{\sqrt{1 - \delta_\lambda^2}}{\pi \sigma f(y)} \left(\frac{\eta^2}{\nu(1 - \delta_\lambda^2)} + 1\right)^{-\left(\frac{\nu}{2}+1\right)},$$

$$\boldsymbol{iii}) \; E(\gamma^2\tau|y) = \delta_\lambda^2 \frac{(y-\xi)^2}{\sigma^2} E(\tau|y) + (1 - \delta_\lambda^2) + \frac{\delta_\lambda(y - \xi)\sqrt{1 - \delta_\lambda^2}}{\pi \sigma^2 f(y)} \left(\frac{\eta^2}{\nu(1 - \delta_\lambda^2)} + 1\right)^{-\left(\frac{\nu}{2}+1\right)},$$

$$\boldsymbol{iv}) \; E(\log \tau|y) = DG\left(\frac{\nu+1}{2}\right) - \log\left(\frac{\eta^2 + \nu}{2}\right) + \left(\frac{\nu+1}{\eta^2 + \nu}\right) \left(\frac{T_{\nu+3}\left(\lambda\eta\sqrt{\frac{\nu+3}{\eta^2+\nu}}\right)}{T_{\nu+1}\left(\lambda\eta\sqrt{\frac{\nu+1}{\eta^2+\nu}}\right)} - 1\right)$$

$$+ \frac{\lambda\eta(\eta^2 - 1)}{\sqrt{(\nu + 1)(\eta^2 + \nu)^3}} \frac{t_{\nu+1}\left(\lambda\eta\sqrt{\frac{\nu+1}{\eta^2+\nu}}\right)}{T_{\nu+1}\left(\lambda\eta\sqrt{\frac{\nu+1}{\eta^2+\nu}}\right)} + \frac{1}{T_{\nu+1}\left(\lambda\eta\sqrt{\frac{\nu+1}{\eta^2+\nu}}\right)} \int_{-\infty}^{M} g_\nu(x) t_{\nu+1}(x) dx.$$

These conditional expectations will be used in EM algorithm given in Section 3.